\documentclass[12pt]{article}
\input epsf.tex
\usepackage{latexsym}
\usepackage{amsmath,amsthm,amsfonts,amssymb}
\usepackage{epsfig}

\newtheorem{pro}{Proposition}[section]
 \newtheorem{thm}[pro]{Theorem}
 \newtheorem{lem}[pro]{Lemma}

 \newtheorem{cor}[pro]{Corollary}
\def\A{{\cal A}}
 \def\B{{\cal B}}
 
\def\E{{\cal E}}
 \def\H{{\mathbb H}}
\def\homothety{{\cal H}}
 \def\G{{\cal G}}
 \def\P{{\cal P}}
 \def\O{{\cal O}}
 
 \def\F{{\cal F}}
 
 \def\T{{\cal T}}
 
 \def\R{{\mathbb R}}
 \def\chix{{\raise.5ex\hbox{$\chi$}}}
 
\def\Z{{\mathbb Z}}

\begin{document}
\title{Periodicity and Circle Packings of the Hyperbolic Plane}
\author{Lewis Bowen\footnote{Research supported in part by NSF Vigre Grant No. DMS-0135345} }

\maketitle
\begin{abstract}
We prove that given a fixed radius $r$, the set of isometry-invariant probability measures supported on ``periodic'' radius $r$-circle packings of the hyperbolic plane is dense in the space of all isometry-invariant probability measures on the space of radius $r$-circle packings. By a periodic packing, we mean one with cofinite symmetry group. As a corollary, we prove the maximum density achieved by isometry-invariant probability measures on a space of radius $r$-packings of the hyperbolic plane is the supremum of densities of periodic packings. We also show that the maximum density function varies continuously with radius.
\end{abstract}

\section{Introduction}
The theory of sphere packings in the Euclidean space $\R^n$ is well-developed and has deep connections to many areas of mathematics. The central question of that theory is: what is the densest packing of Euclidean space by balls of a fixed radius? A sphere or circle packing $P$ is a collection of nonoverlapping congruent balls and the density is defined to be the limit of the relative density of the packing within a ball centered at the origin as the radius of the ball tends to infinity. In symbols,
\begin{equation}
density(P) = \lim_{R \to \infty} density(P \hbox{ in }B_R)
\end{equation}
where $B_R$ denotes the ball of radius $R$ centered at the origin. Of course, the above limit may not exist in which case the density of $P$ is not defined. The terms ``sphere packing'' and ``circle packing''  are unfortunate because they refer to packings by balls, not spheres or circles but they are standard terms so we will use them. 

The central question is considered well-posed because in every dimension $n$ there exists a sphere packing $P$ whose density is greater than or equal to the density of any other sphere packing. Also, the radius is inconsequential; that is the maximum density attained by a sphere packing by balls of radius $r$ is equal to the maximum density attained by a sphere packing by balls of radius $R$ for any $r$ and $R$. This is because homothetic expansion in Euclidean space changes the radius of a packing but not its density.
Also, the limit in the definition of density is independent of the choice of origin.

Except in dimensions 0 and 1 where the problem is trivial, proving a sphere packing $P$ of $\R^n$ is a densest packing is a very hard problem. Fejes Toth [Fe5] gave the first rigorous proof that the hexagonal lattice packing is a densest packing in dimension 2. This is the packing whose set of disk centers coincides with the hexagonal lattice. In 1998, Thomas Hales [Hal] produced a proof that the cannonball packing in dimension 3 is a densest packing. But his proof is so computationally intensive that it has not yet been fully refereed. Beyond dimension 3, no sphere packing has been proven to be densest.

Most sphere packing research has focussed on periodic packings, i.e. packings whose symmetry group is cocompact [CoS]. This is not too restrictive a hypothesis because in any dimension $n$, the supremum over all densities of periodic sphere packings of $\R^n$ is equal to the maximum density possible [BoR2]. However it is not known whether the maximum density in any dimension greater than 3 is attained by a periodic packing. 

Compared with the Euclidean theory, the study of sphere-packings in $n$-dimensional hyperbolic space $\H^n$ is still in its infancy. $\H^n$ is the unique simply connected Riemannian manifold of constant sectional curvature -1. The volume of a ball in $\H^n$ grows exponentially with radius [Rat]. Because of this, the volume of space near the boundary of a large ball is significant relative to the volume of the ball. Consequently, the limit in the definition of density of packing $P$ may depend on the choice of origin [BoR2]. In fact, the limit may exist for some choices of origin but not for others. Also, it is not known whether for every radius $r$ and every dimension $n$ there exists a packing $P$ of $\H^n$ by balls of radius $r$ for which the limit exists independently of origin and such that this common value is maximal over all densities of radius $r$-sphere packings of $\H^n$ (even if we restrict to those packings whose density exists independently of origin). For these and other reasons, some researchers concluded that the sphere-packing problem is not well-posed in hyperbolic space ([Bo1-2], [BoF], [Fe1-5], [FeK]).

One way to avoid this difficulty is to focus on periodic packings. Here we use the term periodic to mean ``having cofinite symmetry'' rather than the more restrictive ``having cocompact symmetry''; these two definitions are equivalent in Euclidean space but different in hyperbolic space. In this case, it is known that the density exists independent of origin and is equal to the relative density of the packing in any fundamental domain for the symmetry group of the packing [FeK]. 

There are no homotheties in hyperbolic space. So, if $R > r >0$, the study of radius $R$-sphere packings of $\H^n$ may be completely different from the study of radius $r$-packings of $\H^n$. It was shown in [BoR1] that in any dimension $n$, for all but a countable set of radii $r$, if $P$ is any periodic radius $r$-sphere packing then there is another another periodic radius $r$-sphere packing $P'$ with greater density than $P$. So, for most radii no periodic packing is ``densest''. 

The difficulties experienced by earlier researchers have mostly been overcome by a new ergodic-theoretic framework for studying packing problems [BoR1]. Oded Schramm suggested an approach to density through the mass-transport principle [BeS]. The authors of [BoR1] combined this idea with new ergodic theory results of Nevo and Stein [Nev,NeS] to produce the new framework. The main idea is to study isometry-invariant probability measures on the space $\P^n_r$ of radius $r$-sphere packings of $\H^n$ in which this space is given the topology of uniform convergence on compact subsets. Such a measure $\mu$ is called ergodic if it is not a nontrivial convex sum of other isometry-invariant probability measures on $\P^n_r$. In this case new ergodic theory results (due to Nevo and Stein [Nev,NeS]) and a relatively small additional argument [BoR2] imply that for $\mu$-almost every packing $P$, the density of $P$ exists (i.e. the above limit converges) independently of origin. Moreover the density of $P$ is equal to the density of $\mu$ which is defined to be the $\mu$-measure of the set of packings that cover the origin. The space of isometry-invariant probability measures $M^n_r$ on $\P^n_r$ is compact under the weak* topology and the density function is continuous on $M^n_r$. So, there exists an ``optimally-dense'' measure. That is, a measure $\mu \in M^n_r$ whose density is at least as great as the density of any other measure in $M^n_r$ and which is ergodic. The density of such a measure $\mu$ is called the optimal density of the sphere of radius $r$ in dimension $n$ and will be denoted here by $D^n(r)$. We call $D^n$ the optimum density function.

The new framework generalizes the study of periodic sphere-packings. This is because, as detailed in the next section, for every periodic radius $r$-sphere packing $P$ there is a unique invariant measure $\mu \in M^n_r$ such that $\mu$-almost every packing $P' \in \P^n_r$ is equal to $P$ up to rigid motion. In this case, $\mu$ is said to be periodic and its density is equal to the density of $P$ [BoR1]. We will also call $\mu$ periodic if it is a convex sum of periodic measures.

Besides proving a framework for the study of hyperbolic sphere-packings, this new machinery can be applied to the more general study of packings of Euclidean (or hyperbolic) space by bodies other than spheres. Indeed the completely saturated-conjecture was recently solved [Bow] using this new technology.

Now that many of the early difficulties in hyperbolic sphere-packing research have been resolved, deeper questions can be asked. This paper focusses on two such questions:

\begin{enumerate}
\item Is the set of periodic measures dense in the space $M^n_r$ of invariant measures on the space $\P^n_r$ of radius $r$-sphere packings?
\item Is the optimum density function $D^n: (0,\infty) \to [0,1]$ continuous?
\end{enumerate}

An affirmative answer to the first question would imply that optimum density is the supremum over all densities of periodic packings. In the Euclidean setting, it is almost trivial that both questions have ``yes'' answers. This is a consequence of two facts: the surface area to volume ratio of a Euclidean cube tends to zero as the size of the cube tends to infinity and the cube is a fundamental domain for a cocompact subgroup of Euclidean isometries (namely $\lambda \Z_n$ for some $\lambda>0$). So any ergodic measure $\mu$ can be approximated in the following way. For a given (large) cube choose $\mu$-almost any packing $P$ and remove from $P$ all balls not contained in the cube. Then identify opposite sides of the cube to obtain a packing of an $n$-torus. Pull this packing back to a periodic packing on Euclidean space by a universal covering map. There is a unique measure supported on the orbit of this periodic packing. This sequence of periodic measures converges to the original measure as the size of the cube tends to infinity.  To approximate an arbitrary invariant measure, represent it as a convex sum (or integral) of ergodic measures and then approximate the ergodic measures. The second question above is trivial in Euclidean space since the optimum density function is constant for every dimension (and equal to the maximum density) because of the existence of homotheties.

But in hyperbolic space, there are no homotheties and the surface area to volume ratio of any large-volume body tends to one instead of zero as the volume tends to infinity. So it is not possible to apply the above methods to the hyperbolic setting. Despite this, we give affirmative answers (theorems 3.1 and 3.2) to both questions in dimension 2.  

Here is an outline of the proofs. We first show that any measure $\mu \in M^2_r$ can be approximated by measures supported on radius $r$-sphere packings that have ``wiggle room''; meaning that the distance between any two sphere centers in such a packing is at least $R$ for some $R>r$ (independent of the spheres). In Euclidean space, this could be accomplished by composing homothetic expansion with resizing the radius of each ball. In hyperbolic space we use a modification of the inverse image of a covering map (of $\H^2$ to itself) branched over the vertices of a regular tesselation. We call such a map a ``branched homothety''; it is introduced and studied in section 6. 

It is almost trivial to show that $D^n$ is left continuous (in every dimension) because for any optimally dense radius $R$-sphere packing measure $\mu$ there is an obvious radius $r < R$ sphere packing measure $\mu_r$ defined by shrinking all the radii of all the packings in the support of $\mu$ from $R$ to $r$. The proof of right-continuity in dimension 2 uses branched homotheties to construct from a radius $r$-sphere packing measure $\mu_r$ a new radius $R$-sphere packing measure $\mu_R$ with $R > r$ whose density is close to the density of $\mu_r$. This answers the second question above (the details are at the end of section 3).

A feature specific to dimension 2 is the existence of cofinite free groups that act on $\H^2$. Fixing such a group $F$, we formulate a question analogous to the first one for ``colorings'' (instead of packings) on $F$ (which we may think of as a geometric space by fixing a Cayley graph for $F$). By a coloring, we mean a function from $F$ to a fixed finite set $K$. We solve this question by explicitly constructing the approximating ``periodic measures'' in section 4. The proof heavily depends on the fact that $F$ has a simply-connected Cayley graph.

If a probability measure $\mu$ on the space of packings is not invariant under all isometries but is invariant under $F$ then we can average it over a fundamental domain for the left action of $F$ on the group of all isometries. This yields a new measure $\A(\mu)$ which is invariant under all isometries. We study this averaging function in section 5. 

To answer the first question, we let $\mu$ be an invariant measure in $M^2_r$. Using branched homotheties, we approximate $\mu$ by a new measure $\mu'$ supported on packings with wiggle room. We use this wiggle room to construct a probability measure $\mu''$ that is supported on packings $P$ whose intersection with a fixed fundamental domain for $F$ occurs in only finitely many configurations. This new measure is $F$-invariant but not fully invariant. We consider the packings in the support of $\mu''$ as colorings of $F$ and approximate $\mu''$ by $F$-invariant probability measures $\mu'''$ supported on packings that have cofinite symmetry. We then average $\mu'''$ over a fundamental domain for the left-action of $F$ on the group of isometries to obtain $\A(\mu''')$. We show that $\A(\mu''')$ is periodic and approximates $\mu$. This answers the first question above (the details are in section 3).

\section{The Framework}

The following background material appears in greater detail in [BoR1-2]. A sphere packing (by spheres of radius $r$ with $0< r < \infty$) is a collection of radius $r$ balls in $\H^n$ so that distinct balls have disjoint interiors. We say that a packing $P$ has radius $r$ if it is radius $r$-sphere packing.

We let $\P=\P^n$ be the space of all locally finite collections $P$ of spheres in $\H^n$ with the topology of uniform convergence on compact subsets. Note that spheres of a collection $P$ in $\P$ are allowed to overlap.

$\P^n$ is a locally compact metrizable space. We let $\P^n_r=\P_r$ be the subspace of all radius $r$ packings. For each $r$, $\P^n_r$ is compact. The group $Isom^+(\H^n)$ of orientation preserving isometries of hyperbolic space acts on $\P^n$ in the natural way. We let $M^n$ be the space of Borel probability measures on $\P^n$ that are invariant with respect to the action of $Isom^+(\H^n)$. By invariant we mean that if $\mu \in M^n, E \subset \P^n$ and $g \in Isom^+(\H^n)$ then $\mu(gE)=\mu(E)$. We consider $M^n$ with the weak* topology which is defined as follows. If $\{\lambda_j\}_j \subset M^n$ then $\lambda_j$ converges to  $\mu$ (weak*) if any only if $\lambda_j(f) \to \mu(f)$ for every continuous function $f: \P^n \to \R$. We let $M^n_r \subset M^n$ be the subspace of measures supported on $\P^n_r$. Then $M^n_r$ is compact for each $r$ by the Banach-Alaoglu theorem. When the dimension $n$ is understood, we write $M_r$ for $M^n_r$.

If $P \in \P^n_r$ has symmetry group $Sym(P) < Isom^+(\H^n)$ that is cofinite (in which case we say $P$ is {\bf periodic}) then its orbit $O(P)= Isom^+(\H^n)P \subset \P^n_r$ is naturally homeomorphic with $Isom^+(\H^n)/Sym(P)$ (by a homeomorphism which intertwines the action of $Isom^+(\H^n)$). Therefore the Haar measure on $Isom^+(\H^n)/Sym(P)$ pulls back to an isometry-invariant measure $\mu$ on $O(P)$ which we may normalize so that $\mu(O(P))=1$. We extend $\mu$ to all of $\P^n_r$ by setting $\mu(\P^n_r - O(P))=0$. Now $\mu \in M^n_r$. If a measure $\mu \in M^n_r$ is obtained from a periodic packing $P$ in this way, we say that $\mu$ is a {\bf periodic measure}. We will also say that $\mu \in M^n_r$ is periodic if it in the convex hull of all such measures. (This definition differs slightly from [BoR1-2] in that now we are allowing the symmetry group of $P$ to be cofinite but not necessarily cocompact).

For convenience we let $\O \in \H^n$ be a point we will call the origin. As suggested by Oded Schramm, we define the density of a measure $\mu \in M^n_r$  by

\begin{equation}
density(\mu)=\mu[\{P \in \P_r| \, \textrm{ the origin is in a ball of } P\}]
\end{equation}

It is shown in [BoR1,2] that if $P$ is a periodic packing with corresponding periodic measure $\mu \in M^n_r$ then the density of $\mu$ is equal to the relative density of $P$ in any fundamental domain of $Sym(P)$.

A measure $\mu \in M^n_r$ is said to be {\bf ergodic} if it is not a nontrivial convex sum of two other measures in $M^n_r$ , i.e. if $\mu = t\lambda_1 + (1-t)\lambda_2$ with $t \in (0,1)$ and $\lambda_i \in M^n_r$ (for $i=1,2$) then $\lambda_1=\lambda_2=\mu$. It was shown in [BoR1,2] that if $\mu \in M^n_r$ is ergodic then $\mu$-almost every packing $P$ is such that the orbit $O(P)$ is dense in the support of $\mu$. Given an ergodic measure $\mu \in M^n_r$, almost every packing $P \in \P^n$ is such that for any point $p \in \H^n$

\begin{equation}
density(\mu) = \lim_{R \to \infty} \frac{vol\big(P \cap B_R(p)\big)}{vol\big(B_R(p)\big)}
\end{equation}

\noindent where $B_R(p)$ is the ball of radius $R$ centered at $p$. Here we have abused notation by identifying $P$ with the closed subset equal to the union of all balls in $P$.

The following lemma and proof are essentially contained in [BoR1,2] although they are not explicitly stated.

\begin{lem}
The function $density:\bigcup_r \, M^n_r \to [0,1]$ is continuous.
\end{lem}

\begin{proof}
Let $A_0$ be the set of all packings $P \in \P$ such that the origin is contained in a ball of $P$. Then, by definition $density(\mu) = \mu(A_0)$ for all $\mu \in M^n$. For $\epsilon > 0$, let $A^+_\epsilon$ be the set of all packings $P \in \P$ such that for any $r$ the origin is less than a distance  $r(1+\epsilon)$ from the center of a radius $r$-ball of $P$. Let $A^-_\epsilon$ be the set of all packings $P \in \P$ such that for any $r$ the origin is at most a distance $r(1-\epsilon)$ of the center of a radius $r$-ball of $P$. It is clear that $A^-_\epsilon \subset A_0 \subset A^+_\epsilon$ and that $A^-_\epsilon$ is closed in $\P^n$ and $A^+_\epsilon$ is open in $\P^n$. By Urysohn's lemma, there exists a continuous function $f_\epsilon: \P \to [0,1]$ such that $f_\epsilon(P) = 1$ for all $P \in A^-_\epsilon$ and $f_\epsilon(P)=0$ for all $P \notin A^+_\epsilon$. Therefore $\mu(A^-_\epsilon) \le \mu(f_\epsilon) \le \mu(A^+_\epsilon)$. 

Fix $\mu \in M^n_r$ (for some $r$) and choose $0<\epsilon < 1$.
It is not hard to show that

\begin{eqnarray}
\mu(A^+_\epsilon) &\le& \frac{vol\big(B_{r+r\epsilon}(\O)\big)}{vol\big(B_r(\O)\big)} density(\mu)\\
\mu(A^-_\epsilon) &=& \frac{vol\big(B_{r-r\epsilon}(\O)\big)}{vol\big(B_r(\O)\big)} density(\mu).
\end{eqnarray}

Suppose $\mu_i \in M^n_{r_i}$ and $\mu_i \to \mu$ (as $i\to \infty$). Choose $i$ large enough so that $|r_i -r|<\epsilon$. By weak* convergence, $\mu_i(f_\epsilon) \to \mu(f_\epsilon)$. Hence 

\begin{eqnarray}
\frac{vol(B_{r-r\epsilon})}{vol(B_r)} \limsup_i density(\mu_i) &=&   \limsup_i \frac{vol(B_{r_i-r_i\epsilon})}{vol(B_{r_i})} density(\mu_i)\\
                              &=&  \limsup_i \mu_i(A^-_\epsilon)\\
                              &\le&  \limsup_i \mu_i(f_\epsilon)\\
                              &=& \mu(f_\epsilon)\\
                              &=&  \liminf_i \mu_i(f_\epsilon)\\
                              &\le& \liminf_i \frac{vol(B_{r_i + r_i\epsilon})}{vol(B_{r_i})} density(\mu_i)\\
                              &=& \frac{vol(B_{r +r\epsilon})}{vol(B_r)} \liminf_i density(\mu_i).
\end{eqnarray}

\noindent where $vol(B_r)$ is the volume of the ball of radius $r$. Taking $\epsilon \to 0$ proves the lemma.

\end{proof}

Since $M^n_r$ is compact there exists an ergodic measure $\mu$ (called an {\bf optimally dense measure}) such that $density(\mu) = \max\{density(\lambda) |, \lambda \in M^n_r\} =: D^n(r)$.

\section{Results}

The main results are:

\begin{thm}
For any radius $r \in (0,\infty)$, the space of periodic measures in $\P^2_r$ is dense in $M^2_r$.
\end{thm}

\begin{cor}
 $D^2(r)=\sup \{density(\mu) | \, \mu \in M^2_r$ is periodic $\}$. 
\end{cor}

\begin{thm}
The optimum density function $D^2:(0,\infty) \to [0,1]$ is continuous.
\end{thm}

\noindent {\bf Remark}: Both theorems apply to dimension 2 only. All higher dimensional cases of the above theorems are unknown. Henceforth we will write $D$ for $D^2$, $M_r$ for $M^2_r$, etc. 

To motivate our techniques, consider the following discretized version of the problem of which theorem 3.1 answers a special case. Let $F$ be a discrete group and let $K$ be a topological space. Let $X=X(F,K)$ denote the space of all functions $\phi: F \to K$ with the topology of uniform convergence on finite subsets. If $K$ is compact, $X$ is a compact metrizable space on which $F$ acts by $f\phi(g)=\phi(f^{-1}g)$ for any $f,g \in F$ and $\phi:F \to K$. We let $M=M(F,K)$ be the space of Borel probability measures on $X$ that are invariant under this action (with the weak* topology). We say that $\mu \in M$ is {\bf periodic} if its support is finite. We believe the following problem to be of intrinsic interest:

\noindent {\bf Problem:} Characerize those discrete groups $F$ with the property that the space of periodic measures is dense in $M(F,K)$ where $K$ is a finite set of at least 2 elements.

In section 4 we prove the following:

\begin{thm} 
If $F$ is a finitely generated free group and $K$ is a finite set, then 
\begin{enumerate}
\item the subspace of periodic measures is dense in $M(F,K)$,
\item suppose $\psi_1,..,\psi_m \in X$ and $E \subset F$ is a finite set. If $\mu \in M(F,K)$ is supported on the set $Y:=\{\phi \in X| \, \exists i \textrm{ such that } \, \phi|_E = \psi_i|_E\}$ then there exists periodic measures $\lambda_j \in M(F,K)$ also supported on $Y$ such that $\lambda_j \to \mu$. 
\end{enumerate}

\end{thm}

\noindent {\bf Remark}: We do not know if this theorem holds when $F$ is the fundamental group of a closed surface of genus $\ge 2$ or when $F$ is the fundamental group of the figure 8 knot complement (for example).

It is possible to apply this theorem to prove theorem 3.1 only because of the existence of cofinite free discrete subgroups of $Isom^+(\H^2)$. However, we first need a nice way of separating the circles of a circle packing. For this purpose, in section 6 we prove

\begin{thm}
There exists a natural family $\{\homothety_{s,a}\}$ of maps from $M_r$ to $M_r$ indexed by integers $(s,a)$ such that there exists a regular hyperbolic polygon with $s$ sides and angles equal to $2\pi/a$ such that the following holds.
\begin{enumerate}
\item For every such $(s,a)$ there exists an $R>r$ such that for any $\mu \in M_r$, $\homothety_{s,a}(\mu)$ is supported on the space of all packings $P$ in $\P_r$ such that any two centers of distinct disks in $P$ are at least a distance $2R$ apart.
\item For fixed $s$, $\homothety_{s,a}$ converges to the identity pointwise as $a \to \infty$.
\item  $density\big(\homothety_{s,a}(\mu)\big)$ converges to $density(\mu)$ as $a \to \infty$.
\end{enumerate}
\end{thm}

 

\noindent We now prove theorems 3.1 and 3.2 given theorems 3.3, 3.4 and the results of section 5.

\begin{proof}(of theorem 3.1):
Let $r>0$ and $\epsilon>0$ be given. Let $\mu$ be a measure in $M^2_r$. Since, for fixed $s\ge3$, $\homothety_{s,a}(\mu)$ converges to $\mu$ as $a \to \infty$ it suffices to show that the measures $\homothety_{s,a}(\mu)$ are approximatable by periodic measures. So we may assume that there is an $R > r$ such that $\mu$ is concentrated on the set of all packings by disks of radius $r$ in which the centers of any two disks are at least a distance $2R$ apart. 

Let $F < Isom^+(\H^2)$ be a discrete cofinite finitely generated free group. It is a standard fact that such subgroups exist. For example, $F$ could be a finite index subgroup of $PSL_2(\Z)$ (where we identify $Isom^+(\H^2)$ with $PSL_2(\R)$ in the standard way [Rat]).

Let $\Delta \subset \H^2$ be a convex fundamental domain for $F$. Consider the set $K_\infty$ of all discrete subsets of $\Delta$. We topologize $K_\infty$ with the topology of uniform convergence on compact subsets. Recall that $X(F,K_\infty)$ is the set of all functions from $F$ to $K_\infty$ with the topology of uniform convergence on compact finite subsets. In order to transfer the problem to a similar problem on the free group $F$, we define $\Phi: \P \to X(F, K_\infty)$ by 

\begin{equation}
\Phi(P)(f) = \big(f^{-1}C(P)\big) \cap \Delta.
\end{equation}

\noindent where $C(P)$ denotes the set of centers of disks in $P$. We claim that $\Phi$ is an $F$-equivariant bijection. So let $f_1, f_2 \in F$. Then

\begin{eqnarray}
\Phi(f_1P)(f_2) &=& f_2^{-1}(f_1 C(P)) \cap \Delta\\
                       &=& \Phi(P)(f_1^{-1}f_2)\\
                       &=&[f_1 \Phi(P)](f_2).
\end{eqnarray}

\noindent Hence, $\Phi(f_1P)=f_1\Phi(P)$ for any $f_1 \in F$, i.e. $\Phi$ is $F$-equivariant. To see that $\Phi$ is invertible, define $\Psi: X(F, K_\infty) \to \P$ by

\begin{equation}
\Psi( \phi) = \cup_{f \in F} \, fC^{-1}\big(\phi(f)\big).
\end{equation}
 
\noindent where $C^{-1}\big(\phi(f)\big)$ denotes the collection of radius $r$ disks whose set of centers is exactly $\phi(f)$. We check that $\Psi \circ \Phi$ is the identity:

\begin{eqnarray}
\Psi(\Phi(P)) &=& \cup_{f \in F} \, f C^{-1}\big( \Phi(P)(f) \big)\\
              &=& \cup_{f \in F} \, f C^{-1}\Big(\big(f^{-1}C(P) \cap \Delta \big)\Big)\\
              &=& \cup_{f \in F} \, P \cap f \Delta\\
              &=& P.
\end{eqnarray}

\noindent Similarly, we can check that $\Phi \circ \Psi$ is the identity.


 $\Phi$ is not continuous because points can cross the boundary of $\Delta$ continuously in $\H^2$. However, it easy to see that $\Phi^{-1} (=\Psi)$ is continuous and $\Phi$ is Borel. Since $\Phi$ and $\Phi^{-1}$ are $F$-equivariant, they induce maps $\Phi_*:M \to M(F, K_\infty)$ and $\Phi^{-1}_*:M(F, K_\infty) \to M_F$ where $M_F$ is the space of measures on $\P$ that are $F$-invariant (but not necessarily $Isom^+(\H^2)$-invariant). By the results of section 5, there is a (natural) continuous map $\A: M_F \to M$ that restricts to the identity on $M$. We will show that $\Phi_*(\mu)$ is approximatable by periodic measures $\lambda_i \in M(F, K_\infty) $ such that for each $i$, $\A \Phi^{-1}_*(\lambda_i)$ is a measure on $\P_r$ (not just on $\P$). By continuity, it follows that $\A \Phi^{-1}_*(\lambda_i) \to \Phi^{-1}_*\Phi_*(\mu) = \mu$. Since the support of $\lambda_i$ is finite (by definition of periodic measure) and $\Phi^{-1}$ commutes with $F$, it follows that $\Phi^{-1}_*(\lambda_i)$ is also supported on a finite set (for each $i$). By lemma 5.3 this implies that $\A \Phi^{-1}_*(\lambda_i)$ is a periodic measure (for each $i$). This completes the theorem (given the existence of the sequence $\{\lambda_i\}$).

The existence of $\{\lambda_i\}$ does not follow immediately from theorem 3.3 since $K_\infty$ is infinite. So we approximate $K_\infty$ by suitably chosen finite sets as follows. First fix a basepoint $p\in \Delta$. Define  $K_x = \{k \in K_\infty | \, k \subset B_x(p)\}$ (where $B_x(p)$ is the closed ball of radius $x$ centered at $p$).


Let $x >0$ and let $\delta>0$ be such that $R-r> \delta$. Since $B_x(p) \cap \Delta$ is compact, there exists a finite partition $\{O_i\}$ of $K_x$ such that the diameter of $O_i$ (with respect to the Hausdorff distance on discrete subsets of $\Delta \cap B_x(p)$) is at most $\delta$ for all $i$. Recall that the Hausdorff distance between two closed sets $E_1$ and $E_2$ is the smallest number $\rho$ such that $E_1$ is contained in the $\rho$-neighborhood of $E_2$ and $E_2$ is contained in the $\rho$-neighborhood of $E_1$. By a partition, we mean that $O_i \cap O_j = \emptyset$ if $i \ne j$ and $\cup_i O_i = K_x$.  Let $k_i \in O_i$ (for each $i$).  Let $\pi_{x, \delta}=\pi: K_\infty \to K_x$ be defined by $\pi(k)=k_i$ if $k \cap B_x(p) \in O_i$. Note that if $k \in K_\infty$ and $x_i >0, \delta_i >0$ are such that $x_i \to \infty$ and $\delta_i \to 0$ as $i \to \infty$ then $\pi_{x, \delta}(k) \to k$ (in the topology of uniform convergence on compact subsets). Let $K(x,\delta) = \{k_1, k_2...\}$. 

The map $\pi_{x, \delta}$ induces maps $L(x,\delta)$ from $X(F, K_\infty)$ to $X(F, K(x,\delta))$ and $L(x,\delta)_*$ from $M(F, K_\infty)$ to $M(F, K(x,\delta))$ defined in the natural way as follows. 

\begin{eqnarray}
L(x,\delta)(\phi)(f) &=& \pi_{x,\delta}\big(\phi(f)\big)\\
L(x,\delta)_*(\lambda)(E) &=& \lambda\big( L(x,\delta)^{-1}(E) \big).
\end{eqnarray}

\noindent for any $\phi \in X\big(F, K(x,\delta)\big)$ and $E \subset X\big(F,K(x,\delta)\big)$.

 Let $\mu_{x, \delta}$ be the image of $\Phi_*(\mu)$ under $L(x, \delta)_*$. Since for any $k \in K_\infty$, $\pi_{x,\delta}(k)\to k$ (as $x \to \infty$ and $\delta \to 0$ simultaneously) $\mu_{x, \delta}$ converges (weak*) to $\Phi_*(\mu)$. Since $\delta < R-r$ and $\mu$ is supported on the subset of $\P_r$ consisting of those packings $P$ whose set of centers of $C(P)$ is also the set of centers of a packing of radius $R$ it follows that $\Phi^{-1}_*(\mu_{x,\delta})$ is supported on $\P_r$ (not just $\P$). To be pedantic, the measure $\Phi^{-1}_*(\mu_{x,\delta})$ is supported on packings $P$ whose set of centers $C(P)$ is also the set of centers of a packing of radius $R - \delta > r$.

Since $\mu_{x, \delta} \in M(F, K(x, \delta))$ and $K(x, \delta)$ is a finite set, it follows from theorem 3.3 that $\mu_{x,\delta}$ is approximatable by periodic measures. Since $\Phi^{-1}_*(\mu_{x,\delta})$ is supported on $\P_r$ (and since this is a local property) the last sentence of theorem 3.3 implies that $\lambda_i$ can be chosen so that $\Phi^{-1}_*(\lambda_i)$ is also supported on $\P_r$ (for all $i$). This finishes the proof of the existence of the sequence $\{\lambda_i\}$ and the theorem.

\end{proof}

{\bf Remark}: The proof above shows that more is true. For any cofinite free group $F$, the supremum of densities of periodic packings (by circles of radius $r$) in which the symmetry groups of the packings are (finite-index) subgroups of $F$ is equal to the optimum density $D(r)$. Therefore, one could approximate $D(r)$ by considering packings in, say, surfaces of the form $\H^2/G$ where $G$ is a finite index subgroup of $PSL_2(\Z)$.

Before beginning the proof of theorem 3.2 we need the following:

\begin{lem}
For every dimension $n$, the function $D^n:(0,\infty) \to [0,1]$ is continuous from the left.
\end{lem}

\begin{proof}
Since $\bigcup_{t \le s \le r} \, M^n_s$ is compact, if $\mu_s$ is an optimally dense measure in $M_s$ (for $t < s< r$), there exists a limit point $\mu_r$ of $\{\mu_s\}$ with $\mu_r$ concentrated on $\P_r$. Since the density function is continuous (lemma 2.1), $\lim_{s \to r} density(\mu_s) = density(\mu_r)$. Since $\mu_s$ is optimally dense, $\lim_{s \to r-} D^n(s) \le D^n(r)$. 

For $s<r$ define the map $\Psi_{r,s}: \P_r \to \P_{s}$ by stating that the center set of $\Psi_{r,s}(P)$ is equal to the center set of $P$. Intuitively, $\Psi_{r,s}$ is a ``shrinking map'' since it fixes the center set but shrinks the radii. Let $\Psi_{r,s*}:M_r \to M_s$ be defined by $\Psi_{r,s*}(\mu)(E)=\mu(\Psi_{r,s}^{-1}(E))$ for any $E \subset \P_{r,s}$ and $\mu \in M_r$. Then, it is easy to see (cf [BoR1]) that

\begin{equation}
density\big(\Psi_{r,s*}(\mu)\big) = \frac{vol(B_s)}{vol(B_r)}density(\mu)
\end{equation}
\noindent where $vol(B_s)$ denotes the volume (in $\H^n$) of the ball of radius $s$. Taking $\mu$ above to be optimally dense, it follows that

\begin{equation}
D^n(s) \ge \frac{vol(B_s)}{vol(B_r)}D^n(r).
\end{equation}

\noindent Since $\lim_{s \to r} \frac{vol(B_s)}{vol(B_r)}=1$ 

\begin{equation}
\lim_{s \to r-} D^n(s) \ge D^n(r).
\end{equation}

Combining this with the first paragraph, we have $\lim_{s \to r-} D^n(s) = D^n(r)$. So $D^n$ is continuous from the left.

\end{proof}

The proof of theorem 3.2 is quite similar to the above lemma although the role of the ``shrinking map'' will now be played by the branched homothety map $\homothety$ (given by theorem 3.4) composed with an ``expanding map''.

\begin{proof}(of theorem 3.2):

Let $r>0$ be given. Because $\bigcup_{r \le R \le R_2 } \, M_R$ is compact if $\mu_R$ is an optimally dense measure in $M_R$ (for all $r< R$), there exists a limit point $\mu_r$ of $\{\mu_R | R > r\}$ such that $\mu_r$ is concentrated on $\P_r$. Since the density function is continuous, $\lim_{R \to r} density(\mu_R) = density(\mu_r)$. Since $\mu_R$ is optimally dense, $\lim_{R \to r+} D^2(R) \le D^2(r)$. 

Let $\epsilon>0$ and $\mu \in M_r$ be given. By theorem 3.4, there exists a map $\homothety: M_r \to M_r$ such that $\homothety(\mu)$ is supported on packings $P \in \P_r$ such that if $p_1, p_2$ are centers of distinct circles in $P$ then $d(p_1,p_2) \ge 2R$ for some $R > r$ and $density\big(\homothety(\mu)\big) \ge (1-\epsilon)density(\mu)$. There is a (continuous) function $\Psi$ (the ``expanding map'') from the support of $\homothety(\mu)$ into $\P_R$ such that the center set of $\Psi(P)$ is the center set of $P$ for all $P$ in the support of $\mu$. Let $B(\mu) \in M_R$ be defined by $B(\mu)(E) = \homothety(\mu)(\Psi^{-1}(E))$ for any $E \subset \P_R$. Then $B(\mu) \in M_R$ and it is not hard to show that

 \begin{equation}
density(B(\mu)) = \frac{area(B_R)}{area(B_r)}density(\homothety(\mu)) \ge \frac{area(B_R)}{area(B_r)}(1 - \epsilon)density(\mu)
\end{equation}
\noindent where $area(B_R)$ denotes the area (in $\H^2$) of a radius $R$-ball. Taking $\mu$ above to be optimally dense, we see that

\begin{equation}
D^2(R) \ge \frac{area(B_R)}{area(B_r)}(1-\epsilon)D^2(r).
\end{equation}
Since $\lim_{R \to r} \frac{vol(B_R)}{vol(B_r)}=1$ and we may choose $\epsilon>0$ arbitrarily close to zero, 

\begin{equation}
\lim_{R \to r+} D^2(R) \ge D^2(r).
\end{equation}

Combining this with the first paragraph, we have $\lim_{R \to r+} D^2(R) = D^2(r)$. So $D^2$ is continuous from the right. Combining this with lemma 3.5 above, $D^2$ is continuous.

\end{proof}

\section{Periodic approximations to random invariant colorings of the free group on $n$ generators}

In this section we prove theorem 3.3.

Let $F_n$ be the free group on $n$ generators. Let $K$ be a finite set with $k$ elements. Let $X=X(F_n,K)$ be the space of functions from $F_n$ to $K$ with the topology of uniform convergence on finite subsets. $X$ is a compact metrizable space. $F_n$ acts on $X$ by $g\phi(h)=\phi(g^{-1}h)$ where $g,h \in F_n$ and $\phi \in X$. Let $M(F_n,K)$ be the set of all Borel probability measures on $X$ that are invariant under the action of $F_n$. We will consider $M(F_n,K)$ with the weak* topology. An element $\mu$ in $M$ is said to be periodic if its support is finite.

It will be convenient to fix a Cayley graph $\Gamma$ for $F_n$. To do this we choose a standard set $\{g_1,..,g_n\}$ of generators of $F_n$. The vertices of $\Gamma$ are the elements of $F_n$ and for every $g \in F_n$ and generator $g_i$ there is an edge between $g$ and $gg_i$ and an edge between $g$ and $gg_i^{-1}$. Note $\Gamma$ is a $2n$-regular tree. The distance $d$ between elements of $F_n$ is defined to be their distance in the graph $\Gamma$ where each edge has length $1$. We let $B_r(g)$ denote the radius $r$-ball centered at $g$ in $\Gamma$ (or in $F_n$).

The topology of $X$ is generated by its cylinder sets which are defined as follows. Let $E$ be a finite subset of $F_n$ and let $\phi$ be any element of $X$. Then $Z(E,\phi):=\{\psi \in X | \, \phi|_E = \psi|_E\}$. We will be particularly interested in the case in which $E = B_r(id)$. We will say that $\lambda \in M(F_n,K)$ is a $(r, \epsilon)$-approximation of $\mu \in M(F_n,K)$ if

\begin{equation}
\Big|\lambda\big(Z(B_r(id), \phi)\big) - \mu\big(Z(B_r(id), \phi)\big)\Big| < \epsilon
\end{equation}

\noindent for all $\phi \in X$. Its obvious that if $\lambda_i$ is a $(r_i,0)$-approximation to $\mu$ and $r_i \to \infty$ as $i \to \infty$ then $\lambda_i$ converges to $\mu$. Also, if $\lambda_i$ is a $(r, \epsilon_i)$-approximation to $\mu$ and $\epsilon_i \to 0$ as $i\to \infty$ then any limit point of $\lambda_i$ is a $(r,0)$-approximation to $\mu$. So for the first half of the theorem it suffices to show that given {\it any} $(r, \epsilon)$ with $r, \epsilon>0$ and any measure $\mu \in M(F_n,K)$ there exists a periodic $(r,\epsilon)$-approximation to $\mu$. Below, we construct such an approximation.


First, we define a directed, labeled finite graph $G$ that captures the radius $r$ information of $X$ (it will soon be clear what this means). The vertices of $G$ are equivalences classes of elements of $X$ where we say that $\phi \sim \psi$ if $\phi$ restricted to $B_r(id)$ is equal to $\psi$ restricted to $B_r(id)$. We denote the equivalence class of $\phi$ by $[\phi]$. There is a directed edge (denoted $e(\phi,i)$ from $[\phi]$ to $[g_i\phi]$ labeled $i$ for each $\phi$ and for each $i$. Note that $[g_i \phi]$ is not necessarily equal to $g_i[\phi]$.

Next, we define a weight function $w$ on $G$. We let $w([\phi]) = \mu\Big(Z\big(B_r(id), \phi)\big)\Big)$ for $\phi \in X$. We let $w\big(e(\phi,i)\big) = \mu\Big(Z\big(B_r(id) \cup B_r(g_i)\big), \phi\Big)$. There are two basic properties satisfied by $w$. First, the sum of the vertex weights is one (because $\mu$ is a probability measure and the corresponding sets partition $X$). Second, for any vertex $v$ and any $i$ ($1 \le i \le n$), the sum of the weights of edges labeled $i$ going out of $v$ is equal to the weight of $v$ (because $\mu$ is $F_n$-invariant). Similarly, for any vertex $v$ and any $i$ ($1 \le i \le n$), the sum of the weights of edges labeled $i$ going {\it into} $v$ is equal to the weight of $v$ These two properties are two sets of linear equations in the weights. If $\mu$ is periodic and ergodic then all of the weights are rational. To get an approximation to $\mu$, we approximate $w$ by a rational function as follows. Because the linear equations defined by the two basic properties have rational coefficients, there is a function $w_p$ on $G$ satisfying both properies and such that $w_p$ takes values in the nonnegative rational numbers and such that $|w_p(v)-w(v)| < \epsilon$ for all vertices $v$. We may also assume that $w_p(e)=0$ for any edge in which $w(e)=0$.

We use the weight function $w_p$ to construct a disjoint union $Q = \cup Q_i$ of quotient graphs $Q_i$ of $\Gamma$ and functions $\phi_i$ from the vertex set of $Q_i$ to $K$. After this is done, we will lift the functions $\phi_i$ to the vertices of $\Gamma$ and thereby obtain functions ${\tilde \phi_i}:V(\Gamma) \to K$ that have finite orbits under the action of $F_n$. Because of these orbits are finite, there is naturally associated to ${\tilde \phi_i}$ a periodic measure $\lambda_i \in M$ (with support equal to the orbit of ${\tilde \phi_i}$). We will show that an appropriate convex sum $\lambda$ of the $\lambda_i$ is an $(r, \epsilon)$-approximation to $\mu$.

Let $N$ be a natural number such that $Nw_p$ takes values only in the nonnegative integers. We will construct a (directed and labeled) graph $Q$ such that each vertex is labeled $[\phi]$ for some $\phi \in X$, each edge is labeled $e(\phi,i)$ for some $\phi \in X$ and $1 \le i \le n$ and each vertex has exactly one edge labeled $i$ going out of it and one edge labeled $i$ going in to it. Also we choose $Q$ so that the number of vertices of $Q$ labeled $[\phi]$ is exactly $Nw_p(\phi)$ and the number of edges of $Q$ labeled $e(\phi,i)$ is exactly $Nw_p(e(\phi,i))$.

Let $Q_0$ be the edge-less labeled graph with $N$ vertices such that exactly $Nw_p([\phi])$ vertices are labeled $[\phi]$. For each $e(\phi,i)$ let $b_1(\phi,i)$ be a set of vertices of $Q_0$ labeled $[\phi]$ that contains exactly $Nw_p(e(\phi,i))$ elements. Also, choose these sets so that if $[g_i\phi] \ne [g_i\psi]$ then $b_1(\phi,i)$ is disjoint from $b_1(\psi, i)$. We can do this by the second basic property. Next we let $b_2(\phi,i)$ be a set of vertices of $Q_0$ labeled $[g_i\phi]$ that contains exactly $Nw_p(e(\phi,i))$ elements. Likewise, we can do this so that if $[\phi] \ne [\psi]$ then $b_2(\phi,i)$ is disjoint from $b_2(\psi,i)$. Next we let $B(\phi,i): b_1(\phi,i) \to b_2(\phi,i)$ be a bijection. Now we let $Q$ be the graph obtained from $Q_0$ by adjoining an edge from $v$ to $B(\phi,i)(v)$ labeled $e(\phi,i)$ for each $v, \phi, i$. It is easy to see that $Q$ satisfies the properties stated above.

$Q$ has finitely many connected components $Q_1,..,Q_m$. For each $Q_i$, let $\phi_i$ be the map from the vertex set of $Q_i$ to $K$ given by $\phi_i(v) = \phi(id)$ where $v$ is labeled $[\phi]$ for some $\phi \in X$. Up to this point we have not used anything about $F_n$ other than it is a finitely generated group but now we will. Because $\Gamma$ is simply connected, there exists a covering map $C_i: \Gamma \to Q_i$. Moreover, we may choose $C_i$ so that it respects the labels. In other words, the edge connecting $C_i(g)$ to $C_i(gg_j)$ is labeled $j$ for any $1 \le j \le n$. Note that $C_i$ is uniquely determined by $C_i(id)$. We let ${\tilde \phi_i} \in X$ be defined by ${\tilde \phi_i}(g) = \phi_i C_i(g)$. It is clear than ${\tilde \phi_i}$ has finite $F_n$-orbit since $Q_i$ is finite. Let $\lambda_i$ be the unique measure in $M(F_n,K)$ whose support is equal to the orbit of ${\tilde \phi_i}$. We let $\lambda = (1/N) \Sigma_i |V(Q_i)|\lambda_i$ where $|V(Q_i)|$ is the number of vertices in the graph $Q_i$. Recall that $N$ is the total number of vertices of $Q$ so $\lambda$ is a probability measure. It is clear that $\lambda$ has finite support and therefore is periodic.

We claim that $\lambda$ is a $(r, \epsilon)$-approximation to $\mu$. Let $N_i([\phi])$ be the number of vertices of $Q_i$ labeled $[\phi]$. Then,

\begin{eqnarray}
\lambda\big(Z(B_r(id),\phi)\big) &=& (1/N) \, \Sigma_i \, |V(Q_i)| \lambda_i \big(Z(B_r(id), \phi)\big)\\
                         &=& (1/N) \, \Sigma_i \, N_i([\phi])\\
                         &=& (1/N)\, Nw_p([\phi])\\
                         &=& w_p([\phi]).
\end{eqnarray}
 
The second equation holds because $\frac{N_i([\phi])}{|V(Q_i)|}$ is the fraction of vertices of $Q_i$ that are labeled $[\phi]$. The third equation holds because $\Sigma_i \, N_i([\phi])$ is the total number of vertices in $\bigcup_i Q_i =Q$ labeled $[\phi]$. By construction this is equal to $Nw_p([\phi])$. Since $\mu(Z\big(B_r(id),\phi)\big) = w([\phi])$ and $\big|w([\phi])-w_p([\phi])\big|<\epsilon$ we are done with the first half. Let $E \subset F$ be a finite set. For $\epsilon>0$ and $r >0$ large enough so that $E \subset B_r(id)$, the construction above yields a periodic measure $\lambda$ that is a $(r,\epsilon)$-approximation of $\mu$ and such that $\lambda\big(Z(B_r(id),\phi)\big)=0$ if $\mu\big(Z(B_r(id),\phi)\big)=0$ (for any $\phi \in X$). This holds because of the requirement that $w_p(e)=0$ whenever $w(e)=0$. This implies the last half of theorem 3.1. 


\section{Averaging partially invariant measures}

In the section we study the properties of a certain averaging process for partially invariant measures. These are necessary results to define the branched homothety $\homothety$ of theorem 3.4. Let $(X,\mu)$ be a probability space on which a locally compact unimodular group $\G$ acts (for the applications of this paper $\G=Isom^+(\H^2)$ and $X \subset \P_r$ is a compact space of packings). Assume that $\mu$ is invariant under a cofinite subgroup $G < \G$. By invariant we mean that $\mu(gE)=\mu(E)$ for all $g\in G, E \subset X$. We let $A(\mu)$ be the probability measure on $X$ defined by

\begin{equation}
\A(\mu)(E) = \frac{1}{vol(\F)} \int_\F \,  \mu(gE)\, dvol
\end{equation}

\noindent where $\F$ is a fundamental domain for the left action of $G$ on $\G$, $vol(\cdot)$ is a Haar measure on $\G$ and $E \subset X$. We will prove that $\A(\mu)$ does not depend on the choice of fundamental domain $\F$, that it is $\G$-invariant, that $\A(\cdot)$ is continuous and the support of $\A(\mu)$ is the closure of the $\G$-orbit of the support of $\mu$.

\begin{lem}
$\A(\mu)$ does not depend on the choice of fundamental domain $\F$.
\end{lem}

\begin{proof}
Let $\F_1, \F_2$ be two fundamental domains for $G$ in $\G$. Let $E \subset X$. Note that $vol(\F_1)=vol(\F_2)$. The lemma follows from the next calculation.

\begin{eqnarray}
\frac{1}{vol(\F_1)} \int_{\F_1} \, \mu(gE) dvol(g) &=& \frac{1}{vol(\F_2)} \Sigma_{h \in G} \, \int_{h\F_2 \cap \F_1} \, \mu(gE) dvol(g)\\
                                                   &=& \frac{1}{vol(\F_2)} \Sigma_{h \in G} \, \int_{h\F_2 \cap \F_1} \, \mu(hgE) \, dvol(g)\\
                                                   &=& \frac{1}{vol(\F_2)} \Sigma_{h \in G} \, \int_{\F_2 \cap h^{-1}\F_1} \, \mu(gE) \, dvol(g)\\
                                                   &=& \frac{1}{vol(\F_2)} \int_{\F_2} \mu(gE) \, dvol(g).
\end{eqnarray}

\noindent The second equation holds from the $G$-invariance of $\mu$.

\end{proof}

 The next lemma is a more general form of lemma 4.2 in [Bow]. The proof is similar.

\begin{lem}
$\A(\mu)$ is invariant under all of $\G$.
\end{lem}

\begin{proof}
Let $E \subset X$ and let $h \in \G$. We must show that $\A(\mu)(hE)=\A(\mu)(E)$. We claim that for some $k$, there exists elements $g_i \in G$ and Borel subsets $\F_i \subset \F$ such that
\begin{enumerate}
\item $\cup_i \, g_i \F_i = \F h$;
\item $vol(\F_i\cap \F_j) = 0$ whenever $i \ne j$;
\item $vol(\cup_i \F_i) = vol(\F)$.
\end{enumerate}

Since $\F$ is a fundamental domain for the (left) action of $G$ on $\G$ there are elements $g_i \in G$ such that $\F h \subset \cup_i g_i \F$. We let $\F_i = g_i^{-1}[g_i \F \cap \F h]$. By definition then, the first part of the claim is true.

Suppose $g \in \F_i \cap \F_j$ for some $i \ne j$. Then $g_i g$ and $g_j g$ are in $\F h \cap G g$. So $g_i gh^{-1}$ and $g_j gh^{-1}$ are in $\F \cap G gh^{-1}$. Since $G g_i gh^{-1} = G g_j gh^{-1}$ and $\F$ is a fundamental domain for $G$, we must have that $g_i gh^{-1}$ and $g_j gh^{-1}$ are on the boundary of $\F$. Hence $g \in g^{-1}_i(\partial\F) h$. Since $g$ is arbitrary, $\F_i \cap \F_j \subset g^{-1}_i (\partial\F) h$. Since $vol(\partial\F) = 0$ the second part of the claim follows. The third part of the claim follows directly from the first two parts and the fact that $vol(\cdot)$ is right-invariant as well as left-invariant since $\G$ is unimodular by hypothesis.

To finish the proof, note
\begin{eqnarray}
\A(\mu)(hE) &=& \frac{1}{vol(\F)} \int_{\F} \mu(ghE) \, dvol(g)\\
          &=&  \frac{1}{vol(\F)} \int_{\cup_i g_i \F_i h^{-1}} \mu(ghE) \, dvol(g)\\
          &=&  \frac{1}{vol(\F)} \sum_i \int_{g_i \F_i h^{-1}} \mu(ghE) \, dvol(g)\\
          &=&  \frac{1}{vol(\F)} \sum_i \int_{\F_i} \mu(g_igE) \, dvol(g)\\
          &=&  \frac{1}{vol(\F)} \sum_i \int_{\F_i} \mu(gE) \, dvol(g)\\
          &=& \frac{1}{vol(\F)} \int_{\F} \mu(gE) \, dvol(g)\\
          &=& \A(\mu)(E).
\end{eqnarray}

The fourth equation holds by the change of variables theorem and because $vol(\cdot)$ is right-invariant. The fifth equation holds because $\mu$ is invariant under $G$. The sixth equation holds because $\F_i$ partitions $\F$ (off of a set of measure zero).
\end{proof}

The next lemma is clear from the definitions.

\begin{lem} The support of $\A(\mu)$ is the closure of the $\G$-orbit of $\mu$. Therefore if $X \subset \P_r, \G = Isom(\H^n)$ and every element in the support of $\mu$ has a symmetry group of finite index in $G$ then every element in the support of $\A(\mu)$ has a cofinite symmetry group in $\G$ so $\A(\mu)$ is periodic.
\end{lem}

\begin{lem} $\A$ is continuous as a function from $M_G$ (the space of $G$-invariant probability measures on $X$) to $M$ the space of all probability measures on $X$. Here both spaces are considered with the weak* topology.
\end{lem}

\begin{proof}
Suppose that $\{\mu_i\}$ is a sequence of measures in $M_G$ that converges to $\mu \in M_G$. Let $f:X \to [0,\infty)$ be a continuous function. We must show that $\A(\mu_i)(f) \to \A(\mu)(f)$. Since $X$ is compact, $f$ has a finite upper bound $M$ say. Let $\epsilon > 0$. Because $\G$ is locally compact there is a compact set $\F_\epsilon \subset \F$ such that $\frac{vol(\F_\epsilon)}{vol(\F)} > 1-\epsilon/M$. Since $\F_\epsilon$ is compact, there is an $I$ such that $i>I$ implies that $\Big|\mu_i\big(g_*(f)\big)-\mu\big(g_*(f)\big)\Big| < \epsilon$ for all $g \in \F_\epsilon$. We compute

\begin{eqnarray}
\Big|\A(\mu)-  \lim_{i \to \infty} \A(\mu_i)(f)\Big| &=& \Big|\frac{1}{vol(\F)} \int_\F \, \mu(g_*f) \, dvol(g)\\
                                           &-&  \lim_{i \to \infty} \frac{1}{vol(\F)} \int_{\F} \, \mu_i(g_*f) \, dvol(g)\Big|\\
                                &\le& \lim_{i \to \infty} \frac{1}{vol(\F)} \int_\F \, |\mu(g_*f) - \mu_i(g_*f)| \, dvol(g)\\
                                &\le& \lim_{i \to \infty} \frac{1}{vol(\F)} \int_{\F_\epsilon} |\mu(g_*f) - \mu_i(g_*f)| \, dvol(g) + \epsilon\\
                                &\le& \frac{1}{vol(\F)} \int_{\F_\epsilon} \epsilon \, dvol(g) + \epsilon \\
                                &\le& 2\epsilon.
\end{eqnarray}

Because $\epsilon$ and $f$ are arbitrary, $\A(\mu_i) \to \A(\mu)$. Since $\mu_i$ and $\mu$ are arbitrary, the result follows.

\end{proof}

\section{Branched Homotheties}

In this section we introduce a tool we call ``branched homotheties'' and prove theorem 3.4. First, we let $ s, a \ge 2$ be integers such that there exists a regular polygon in $\H^2$ with $s$ edges and interior angles equal to $2\pi/a$. Let $\T$ be a regular tiling of the plane by this regular polygon and let $\T_b$ be the regular tiling of the plane by the polygon that has $s$ edges and angles equal to $\pi/a$. There is a branched covering map $\B$ from the plane to itself which takes $\T_b$ onto $\T$ (and is branched over the vertices of $\T$). To define this map, let $\tau_1$ be a tile of $\T_b$ and let $\tau_2$ be a tile of $\T$ and consider the projective disk model of the hyperbolic plane. We may, after moving $T_b$ and $T$ by (possibly different) isometries, assume that the centers of $\tau_1$ and $\tau_2$ are both equal to the origin and that there is a Euclidean similarity that takes $\tau_1$ onto $\tau_2$. To see that this is possible, recall that geodesics in the projective disk model are Euclidean straight lines so both $\tau_1$ and $\tau_2$ are regular Euclidean polygons (with the same number of sides). We let $\B: \tau_1 \to \tau_2$ be the Euclidean similarity that maps $\tau_1$ onto $\tau_2$. Now extend $\B$ to all of $\H^2$ be requiring that for any $\gamma_1$ in the symmetry group of $\T_b$, $\B(\gamma_1 \tau_1)= \gamma_2 \tau_2$ for some $\gamma_2$ in the symmetry group of $\T$. We also require $\B$ to be a local homeomorphism away from the vertices of $T_b$. These properties uniquely define $\B$. The map $\homothety$ will be constructed from the inverse of $\B$.

\begin{lem}
The map $\B$ restricted to $\tau_1$ (a tile of $T_b$) (strictly) decreases distances. Precisely, there exists a constant $c<1$ (depending only on $(s,a)$) such that for $p_1,q_1\in\tau_1$, $d\big(\B(p),\B(q)\big) \le c d(p,q)$.
\end{lem}

\begin{proof}
The element of arclength in the projective disk model is such that the (hyperbolic) length of a vector $v$ based at a point $x$ is
\begin{equation}
\frac{\big[ (1-|x|^2)|v|^2 + <x,v>^2 \big]^{\frac{1}{2}} }{1- |x|^2}
\end{equation}
\noindent where $|x|, |v|, <x,v>$ denotes the Euclidean length of $x$ and $v$ and the Euclidean dot product of $x$ and $v$ respectively ([Rat] theorem 6.1.5). Since $\tau_1$ is compact, it suffices to show that for any $|x|<1, v \in \R^2$ and scalar  $0<k < 1$,

\begin{equation}
\frac{\big[ (1-|kx|^2)|kv|^2 + <kx,kv>^2 \big]^{\frac{1}{2}} }{1- |kx|^2} < \frac{\big[ (1-|x|^2)|v|^2 + <x,v>^2 \big]^{\frac{1}{2}} }{1- |x|^2}.
\end{equation}

\noindent To this end, we fix $x$ (with $|x|<1$) and $v \in \R^2 , v\ne 0$ and define
\begin{equation}
f(k) = \frac{\big[ (1-|kx|^2)|kv|^2 + <kx,kv>^2 \big]^{\frac{1}{2}} }{1- |kx|^2}.
\end{equation}

We need only show that $f$ is strictly monotonic increasing when $1 > k>0$. So we show its first derivative is positive:
\begin{eqnarray}
f'(k) &=& \frac{ (1-|kx|^2)(1/2)\Big( 2k|v|^2 - 4k^3|x|^2|v|^2 + 2k<x,v>^2 \Big) }{ (1-|kx|^2)^2 \Big( (1-|kx|^2)|kv|^2 + <kx,kv>^2 \Big)^{\frac{1}{2}} }\\
      &-& \frac{ \Big( (1-|kx|^2)|kv|^2 + <kx,kv>^2 \Big)(-2k|x|^2) }{(1-|kx|^2)^2 \Big( (1-|kx|^2)|kv|^2 + <kx,kv>^2 \Big)^{\frac{1}{2}} }\\
      &=& \frac{ (1-|kx|^2)\Big( k|v|^2 + k<x,v>^2 \Big) + <kx,kv>^2(2k|x|^2)  }{(1-|kx|^2)^2 \Big( (1-|kx|^2)|kv|^2 + <kx,kv>^2 \Big)^{\frac{1}{2}} }\\
      &>& 0.
\end{eqnarray}
\end{proof}

\begin{cor}
Let $R=r/c > r$. Suppose $\B(p_1)=p_2$ and $\B(q_1)=q_2$ ($p_1$ and $q_1$ are not assumed to be in the same tile of $T_b$). Suppose also that $d(p_2,q_2)\ge 2r$. Then $d(p_1,q_1) \ge 2R$.
\end{cor}

\begin{proof}
This is immediate from the above lemma.
\end{proof}

Fix a radius $r$. For a collection $P$ of circles, let $C(P)$ denote the set of centers of $P$. Conversely, if $P'$ is a discrete set of points, let $C^{-1}(P')$ denote the collection of all radius $r$ circles with centers in $P'$. As $\B$ is a map from $\H^2$ to $\H^2$, we may also consider it (and its inverse) as a map from discrete points to discrete point sets. So consider the map $\homothety_0:=C^{-1}\B^{-1}C: \P \to \P$. If $P \in \P_r$ then $\B^{-1}C(P)$ is not necessarily the set of centers of a packing in $\P_R$ for some $R > r$. There may be two points $p,q \in \B^{-1}C(P)$ such that $d(p,q) < 2R$.  By the corollary above, this can only happen when $\B(p)=\B(q)$. Fixing $p$ for the moment, we see that if $q$ is such that $\B(p)=\B(q)$ but $q\ne p$, then the closest $q$ can be to $p$ (without coinciding with $p$) is when the tile of $\T_b$ containing $q$ intersects the tile of $\T_b$ containing $p$ in a vertex $v$. In this case, the midpoint of the geodesic segment from $p$ to $q$ is $v$ (because $\B(p)=\B(q)$ and $\B$ is branched at $v$ with branch index 2). In particular, if $d(p,q) < 2R$ then $d(p,v) < R$. But this implies that $d(\B(p),\B(v)) < r$.

To compensate for this we define $\E: \P \to \P$ by 

\begin{equation}
\E(P) = C^{-1}\Big( C(P) - \cup_{v \in V(\T)} B_r(v) \Big).
\end{equation}
In other words, $\E(P)$ is the same as $P$ but with all disks whose centers are within a distance $r$ from a vertex of $\T$ erased. Now we have that if $P \in \P_r$ then $\B^{-1}C(\E(P))$ is the center set of a packing in $\P_R$ for some $R > r$ (independent of $P$).  

Clearly $\E$ commutes with the symmetry group of $\T$ and from the definition of $\B$, for every element $g_b$ in the symmetry group of $\T_b$ there is an element $g$ in the symmetry group of $\T$ such that $\B(g_b p)=g\B(p)$ for all points $p \in \H^2$. These two facts combined implies that the map $P \to \homothety_0 \E(P)$ induces a map $\homothety_{0*} \E_* : M \to M(Sym(\T_b))$ where $M(Sym(\T_b))$ is the space of probability measures on $\P$ that are invariant under the symmetry group of $\T_b$. We define the averaging map $\A: M(Sym(\T_b)) \to M$ as in section 5 to get a map $\homothety: M_r \to M_r$ defined by

\begin{equation}
\homothety(\mu) = \A \homothety_{0*} \E_* (\mu). 
\end{equation}

\noindent Note that $\homothety$ depends on $s$ and $a$. When we wish to acknowledge this dependence we will write $\homothety_{s,a} := \homothety$.

We now turn our attention to showing that given $\mu \in M_r$ and $s \ge 3$, $\homothety_{s,a}(\mu)$ converges to $\mu$ as $a \to \infty$. For convenience, we formulate sufficient conditions for convergence in $M_r$. For $\rho, \epsilon >0$ and $P \in \P_r$ we let $N^\rho_\epsilon(P)$ denote the space of all packings $Q \in \P_r$ so that the Haussdorf distance between the center sets of $P$ and $Q$ restricted to the radius $\rho$ ball centered at the origin is less than $\epsilon$. For a fixed $P$, the collection of all neighborhoods $N^\rho_\epsilon(P)$ forms a basis for the neighborhood topology at $P$. We will say that $\lambda \in M_r$ is a $(\rho, \epsilon)$-approximation to $\mu \in M_r$ if for every packing $P \in \P_r$

\begin{equation}
\Big|\lambda \big(N^\rho_\epsilon(P)\big) - \mu\big(N^\rho_\epsilon(P)\big)\Big| < \epsilon.
\end{equation}

\noindent It is not hard to show that if $\lambda_i$ is a $(\rho, \epsilon_i)$-approximation to $\mu \in M_r$ for all $\rho \le \rho_i$ and $\rho_i \to \infty$ and $\epsilon_i \to 0$ as $i \to \infty$ then $\lambda_i \to \mu_i$ in the weak* topology. 

\begin{lem}
Given $\mu \in M_r$ and $s\ge 3$, $\homothety_{s,a}(\mu)$ converges to $\mu$ as $a \to \infty$.
\end{lem}

\begin{proof}
It suffices to show that for every $\rho>0, \epsilon>0$ if $a$ is large enough then $\homothety_{s,a}(\mu)$ is a $(\rho,\epsilon)$-approximation to $\mu$. So let $(\rho, \epsilon)$ be given. For $\delta>0$ define the function $f_\delta$ on $\P_r$ by

\begin{equation}
f_\delta(P) = \mu\big(N^\rho_{\epsilon+\delta}(P)\big) - \mu\big(N^\rho_{\epsilon-\delta}(P)\big).
\end{equation}

\noindent It is not hard to show that $f_\delta$ is continuous on $\P_r$. Since $\P_r$ is compact $\max_P \, f_\delta(P)$ exists and decreases monotonically to zero as $\delta \to 0$. So we choose $\delta>0$ so that $f_\delta < \epsilon$. 


Let $\tau_{s,a}$ denote the regular polygon in $\H^2$ with $s$ sides and angles equal to $2\pi/a$. Let $G_{s,a}$ be the group of symmetries of a regular tiling by $\tau_{s,a}$. Let $\F=\F_{s,a}$ be a fundamental domain for the left action of $G_{s,a}$ on $Isom^+(\H^2)$. We will show that when $a$ is sufficiently large, there is a subset $\F'=\F'_{s,a} \subset \F_{s,a}$ such that the volume of $\F_{s,a} - \F'_{s,a}$ is less than $\epsilon$ and for all $g \in \F'_{s,a}$ and $P \in \P_r$,
\begin{equation}
gN^\rho_{\epsilon - \delta}(P) \subset \homothety_0 \E \big (gN^\rho_\epsilon (P)\big) \subset gN^\rho_{\epsilon + \delta}(P).
\end{equation}

Given this claim, the lemma follows from the following calculation (we have written $N_\epsilon$ for $N^\rho_\epsilon$):

\begin{eqnarray}
&& \Big|\mu\big(N_\epsilon(P)\big) - \homothety(\mu)\big(N_\epsilon(P)\big) \Big|\\
&=& \Big| \mu\big(N_\epsilon(P)\big) - \frac{1}{vol(\F)}\int_\F \mu\Big(\homothety_0 \E\big(gN_\epsilon(P)\big)\Big) dvol(g) \Big|\\
&=&  \frac{1}{vol(\F)}\Big| \int_\F \mu\big(gN_\epsilon(P)\big) - \mu(\homothety_0 \E\big(gN_\epsilon(P)\big) dvol(g) \Big|\\
&\le& \frac{1}{vol(\F)}\int_{\F'} \Big| \mu\big(gN_\epsilon(P)\big) - \mu(\homothety_0 \E\big(gN_\epsilon(P)\big) \Big| dvol(g)  + 1- \frac{vol(\F')}{vol(\F)} \\ 
&\le& \frac{1}{vol(\F)}\int_{\F'}  \mu\big(gN_{\epsilon+\delta}(P)\big) - \mu\big(gN_{\epsilon-\delta}(P)\big)\Big\} \, dvol(g) + \frac{\epsilon}{vol(\F)}\\
&=& \frac{1}{vol(\F)}\int_{\F'}  \mu\big(N_{\epsilon+\delta}(P)\big) - \mu\big(N_{\epsilon-\delta}(P)\big) \, dvol(g) + \frac{\epsilon}{vol(\F)}\\
&=& \frac{vol(\F')f_\delta(P) + \epsilon}{vol(\F)}\\
&\le& \epsilon \frac{vol(\F) - \epsilon + 1}{vol(\F)}\\
&\le& \epsilon\Big(1+ \frac{1}{(s-2 - 2s/a)\pi         }\Big). 
\end{eqnarray}

\noindent Here we have denoted Haar measure on $Isom^+(\H^2)$ by $vol(\cdot)$ and we have normalized it so that the volume of $\F$ is equal to the area of a fundamental domain for $G_{s,a}$ in $\H^2$ (which we computed to be $(s-2 - 2s/a)\pi$). In equation (68) we used that $\mu$ is $Isom^+(\H^2)$-invariant. When $a$ is large the last line above is at most $2\epsilon$ which implies $\homothety_{s,a}(\mu)$ is an $(2\epsilon, \rho)$-approximation to $\mu$. Since $\epsilon, \rho, s$ are arbitrary this implies the lemma.

As before, we will use the projective disk model and assume that the center of $\tau_{s,a}$ (and $\tau_{s,2a}$) coincides with the origin. We let $\B_{s,a}: \H^2 \to \H^2$ be defined by $\B_{s,a}(x)=kx$ where $0<k<1$ is a constant (depending on $s,a$) chosen so that $\B_{s,a}$ maps $\tau_{s,2a}$ onto $\tau_{s,a}$. It is clear that as $a\to \infty$ the $k$ in the definition of $\B_{s,a}$ converges to 1. This implies that $\B_{s,a}$ converges to the identity uniformly on compact subsets. 

We let $\rho_{s,a}>0$ be the smallest radius so that if $p,q$ are points in the polygon $\tau_{s,2a}$ that are not within $\rho_{s,a}$ of a vertex (of $\tau_{s,a}$) then

\begin{equation}
d\big(\B(p),\B(q)\big) \ge d(p,q) - \delta.
\end{equation}

\noindent Because $\B_{s,a}$ converges to the identity, the area of the $\rho_{s,a}$ neighborhood of the vertices of $\tau_{s,2a}$ in $\tau_{s,2a}$ converges to zero as $a \to \infty$. In fact the hyperbolic distance from the origin to a vertex of $\tau_{s,2a}$ minus $\rho_{s,a}$ goes to $\infty$ as $a \to \infty$. Therefore, for all $a$ large enough, the area of the $\rho_{s,a} + \rho + r$ neighborhood of the vertices of $\tau_{s,2a}$ in $\tau_{s,2a}$ is less than $\epsilon$. 


We let $\F'_{s,a}$ be the set of all $g\in \F_{s,a}$ such that $g\O$ (where $\O$ is the origin) is not contained in the $\rho_{s,a} + \rho + r$ neighborhood of any vertex of the tiling by $\tau_{s,a}$. It follows from the above that the volume of $\F'_{s,a}$ is at least equal to the volume of $\F_{s,a} - \epsilon$. Therefore if $g$ is in $\F'_{s,a}$ then

\begin{equation}
\E gN^\rho_\epsilon(P) = gN^\rho_\epsilon(P)
\end{equation}

\noindent since $g\O$ is far enough from a vertex that nothing in the $\rho$ neighborhood of $g\O$ is erased by $\E$. By the above it also follows that the center sets of $\homothety_0(P)$ and $P$ restricted to the $\rho$ neighorhood of $g\O$ are within a distance $\delta$ apart (in the Haussdorf metric). This is because $g\O$ is far enough away from a vertex (specifically it is at least a distance $r + \rho_{s,a}$ away). So,

\begin{equation}
gN^\rho_{\epsilon - \delta}(P) \subset \homothety_0 \E \big (gN^\rho_\epsilon (P)\big) \subset gN^\rho_{\epsilon + \delta}(P).
\end{equation}

\noindent This finishes the claim and the lemma.
\end{proof}

Theorem 3.4 follows immediately from corollary 6.3 and lemma 6.3. 

\section{Questions}

\begin{enumerate}

\item Other bodies: Fix a connected compact set $B$ equal to the closure of its interior. A $B$-packing of $\H^n$ is collection of nonoverlapping congruent copies of $B$. As in this paper, we consider isometry-invariant measures on the space of $B$-packings. We believe that the main results of this paper (theorem 3.1 and 3.2) hold for $B$-packings when $B$ is convex with only minor changes in the proofs. However, when $B$ is nonconvex the branched homothety maps $\H_{s,a}$ are not well-defined. Indeed, the optimum density function is discontinuous on the space of bodies under the Hausdorff topology (even in Euclidean space). So none of the proofs in this paper can be applied to $B$-packings. Is the space of periodic $B$-packing measures dense in the space of all isometry-invariant $B$-packing measures? In Euclidean space, the answer is yes with essentially the same proof as given in the introduction.

\item Periodic approximation: For what other groups $F$ does theorem 3.3 hold? One necessary condition is that the group be residually finite. To see this consider the random coloring $\phi: F \to \{0,1\}$ such that for each element $f\in F$, $\phi(f)$ is an independent random variable that takes the value $1$ with probability 1/2. By definition $F$ is not residually finite if and only if the intersection $L$ of all finite index subgroups of $F$ is not the identity. If this is the case, any periodic measure $\mu \in M(F,K)$ is supported on periodic colorings $\psi$ which must be constant on $L$. So the distribution of $\phi$ is not in the closure of the periodic measures. It appears that stronger subgroup separability properties may also be necessary. Only recently has it been proven that the fundamental group of the figure eight knot complement is subgroup separable on geometrically finite subgroups ([ALR]) so this appears to be a hard but interesting problem.
 
\item Is $D^n$ monotonic? There is a known monotonic increasing upper bound in some dimensions; the density of a radius $r$-packing in any Voronoi region is at most equal to the density of the region occupied by the radius $r$-balls centered at the vertices a regular simplex with sidelength $2r$ in the simplex. This was proven by Fejes Toth [Fe5] in dimension 2 and by Karoly Bor\"ozcky [Bo2] in all dimensions. This bound is known to be monotonic increasing in dimensions 2 (proven by Krammer cited in [Fe5, section 37]) and 3 ([BoF]) and in all sufficiently large dimensions [Mar] (although an explicit lower bound on the dimension has not been provided). In dimension 2, this upper bound is attained for a countable discrete set of radii [Fe5]. Precisely, for every $k \ge 7$ there is a radius $r_k$ such that the equilateral triangle of sidelength $2k$ has interior angles equal to $2\pi/k$. A circle-packing by disks of radius $r_k$ exists such that the set of circle centers coincides with the vertices of an regular tesselation by equilateral triangles of sidelength $2r$. This packing is periodic and the unique invariant measure supported on its orbit is optimally dense [BoR1]. If $D^2$ is not monotonic, it would seem likely that some $r_k$ is a local maxima. 

\item Continuity beyound dimension 2: Is the density function $D^n$ continuous when the dimension $n>2$? The techniques presented in this paper may generalize to dimension 3 where it is known that there are plenty of complete finite volume hyperbolic 3-manifolds branch covering other complete finite volume hyperbolic 3-manifolds. However, when the dimension $n>3$ it is not known whether there exists a complete finite volume hyperbolic $n$-manifolds that nontrivially branch covers another complete finite volume hyperbolic $n$-manifold. Gromov and Thurston ([GrT]) provide some evidence that such branched covers may not exist if the branch locus is totally geodesic. In dimension 4, there is a radius $r_D>0$ and a radius $r_D$-sphere packing $P$ that realizes the Bor\"oczky bound given in question 3 ([Dav]). It is the packing associated to the 120-cell. If $D^4$ is not continuous, a natural place to look for a discontinuity is at $r_D$.

\item Continuity at infinity: There is a natural extention of $D^n$ to $[0,\infty]$. Here a radius $0$-sphere packing of $\H^n$ means a sphere-packing of Euclidean $n$-dimensional space by balls of a fixed radius. A radius $\infty$-packing of $\H^n$ is a horoball packing. Greg Kuperberg pointed out to us that $D^n$ is continuous at $0$ (for any dimension $n$). Here is a sketch of the argument: suppose $P$ is saturated packing of $\H^n$ by balls of radius $\epsilon$. Saturated means every point of $\H^n$ is within a distance $\epsilon$ of some ball of $P$. There is a $k$-bi-Lipschitz homemorphism between $\H^n$ and a piecewise Euclidean manifold so that each piece of the manifold is mapped onto a Voronoi cell [CoS] of $P$ (and so that $k$ tends to one as $\epsilon$ tends to zero). Let $P'$ be a densest sphere packing on the piecewise Euclidean manifold (by balls of a small radius) conditioned on the property that no balls of $P'$ overlaps the boundary of the pieces. Then the set of sphere centers of $P'$ can be pushed-forward to a set $S$ of sphere centers on $\H^n$. A radius $\epsilon'$ can be chosen so that the collection $P''$ of balls centered at $S$ with radius $\epsilon'$ do not overlap. The process of constructing $P''$ from $P$ can be done in a canonical isometry-invariant way. So an invariant radius $\epsilon$-sphere packing measure $\mu$ pushes forward to an invariant radius $\epsilon'$-sphere packing measure $\mu'$. As $\epsilon$ tends to 0, the density of $\mu'$ converges to the maximum density of a Euclidean sphere packing.

Results of Fejes Toth [Fe5] imply that the densest horoball packing in dimension 2 is the packing associated to the modular group. This packing can be obtained from a regular tesselation by ideal triangles by placing the centers of the horoballs at the ideal vertices of the triangles in such a way that all three horoballs centered at the vertices of given triangle are mutually tangent. The density of this packing is $3/\pi$. Because $D^2(r_k)$ (see question 3 for definition of $r_k$) converges to $3/\pi$ as the $k$ tends to infinity it follows that $D^2$ is continuous at infinity as well. If $n>2$ is $D^n$ continuous at infinity? Of course, continuity at infinity is the same as left-continuity at infinity. Lemma 3.6 proves left-continuity at all finite values in all dimensions but the argument does not extend to horoball packings.

\end{enumerate}

{\bf Acknowledgements}: I would like to thank my advisor, Charles Radin for first suggesting this problem and for many helpful conversations. I would also like to thank Oded Schramm for suggesting the approach to density through invariant distributions without which this work could not have gotten started.

Department of Mathematics, One Shields Avenue, Davis, CA, 95616

lbowen@math.ucdavis.edu

\end{document}